\documentclass{amsart}
\usepackage{amssymb}
\usepackage{amsfonts}
\usepackage{amsmath}
\usepackage{graphicx}
\usepackage[pagebackref]{hyperref}
\newtheorem{thm}{Theorem}[section]
\newtheorem{corollary}[thm]{Corollary}
\newtheorem{lemma}[thm]{Lemma}
\newtheorem{proposition}[thm]{Proposition}
\newtheorem{conjecture}[thm]{Conjecture}

\theoremstyle{definition}

\newtheorem{example}[thm]{Example}

\theoremstyle{remark}

\newtheorem{claim}{Claim}

\DeclareMathOperator{\htt}{ht} \DeclareMathOperator{\qf}{qf}

 \DeclareMathOperator{\Max}{Max}
\DeclareMathOperator{\Cl}{Cl} \DeclareMathOperator{\Pic}{Pic}

\def\1{{\rm (1)}}
\def\2{{\rm (2)}}
\def\3{{\rm (3)}}
\def\4{{\rm (4)}}
\def\5{{\rm (5)}}
\def\i{{\rm (i) }}
\def\ii{{\rm (ii) }}
\def\iii{{\rm (iii) }}

\begin{document}
\parindent=.8cm

\begin{center}
{\bf\huge $t$-Class semigroups of integral domains}
\bigskip

{By {\it S. Kabbaj} and   {\it A. Mimouni} at Dhahran}
\bigskip

{\small\tt to Marco Fontana on the occasion of his sixtieth
birthday}
\rule{4cm}{0.1mm}
\bigskip
\end{center}

{\bf Abstract.} The $t$-class semigroup of an integral domain is the semigroup of
the isomorphy classes of the $t$-ideals with the operation induced
by ideal $t$-multiplication. This paper investigates ring-theoretic
properties of an integral domain that reflect reciprocally in the
Clifford or Boolean property of its $t$-class semigroup. Contexts
(including Lipman and Sally-Vasconcelos stability) that
suit best $t$-multiplication are studied in an attempt to generalize
well-known developments on class semigroups. We prove that a Pr\"ufer $v$-multiplication domain (PVMD) is of Krull type (in the sense of Griffin \cite{Gr2}) if and only if its $t$-class semigroup is Clifford. This extends Bazzoni and Salce's results on valuation domains \cite{BS} and Pr\"ufer domains \cite{Ba1,Ba2,Ba3,Ba4}.
\bigskip

\pagestyle{myheadings}
\markboth{\it S. Kabbaj and A. Mimouni, $t$-Class semigroups}
{\it S. Kabbaj and A. Mimouni, $t$-Class semigroups}

\section{\bf Introduction}

The class semigroup of an integral domain $R$, denoted ${\mathcal
S}(R)$, is the semigroup of nonzero fractional ideals modulo its
subsemigroup of nonzero principal ideals \cite{BS,ZZ}. We define the
$t$-class semigroup of $R$, denoted ${\mathcal S}_{t}(R)$, to be the
semigroup of fractional $t$-ideals modulo its subsemigroup of
nonzero principal ideals, that is, the semigroup of the isomorphy
classes of the $t$-ideals of $R$ with the operation induced by
$t$-multiplication. One may regard ${\mathcal S}_{t}(R)$ as the
$t$-analogue of ${\mathcal S}(R)$, exactly, as the class group
$\Cl(R)$ is the $t$-analogue of the Picard group $\Pic(R)$. We have
$\Pic(R)\subseteq \Cl(R)\subseteq{\mathcal S}_{t}(R)\subseteq
{\mathcal S}(R)$. The first and third containments turn into
equality in the class of Pr\"ufer domains as the second does so in
the class of Krull domains. More details on the $t$-operation are
provided in the next section.

A commutative semigroup $S$ is said to be Clifford if every element
$x$ of $S$ is (von Neumann) regular, i.e., there exists $a\in S$
such that $x^{2}a=x$. The importance of a Clifford semigroup $S$
resides in its ability to stand as a disjoint union of subgroups
$G_e$, where $e$ ranges over the set of idempotent elements of $S$,
and $G_e$ is the largest subgroup of $S$ with identity equal to $e$
(Cf. \cite{Ho}). The semigroup $S$ is said to be Boolean if for each
$x\in S$, $x=x^{2}$.

Divisibility properties of $R$ are often reflected in group or
semigroup-theoretic properties of $\Cl(R)$ or ${\mathcal S}(R)$. If
$R$ is a Pr\"ufer domain, $\Cl(R)$ equals its ideal class group, and
then $R$ is a B\'ezout domain if and only if $\Cl(R)=0$. If $R$ is a
Krull domain, $\Cl(R)$ equals its usual divisor class group, and
then $R$ is a UFD  if and only if $\Cl(R)=0$. So an integral domain
$R$ is a UFD if and only if every $t$-ideal of $R$ is principal.
Trivially, Dedekind domains (resp., PIDs) have Clifford (resp.,
Boolean) class semigroup. In 1994, Zanardo and Zannier proved that
all orders in quadratic fields have Clifford class semigroup
\cite{ZZ}. They also showed that the ring of all entire functions in
the complex plane (which is B\'ezout) fails to have this property. In
1996, Bazzoni and Salce investigated the structure of ${\mathcal
S}(V)$ for any arbitrary valuation domain $V,$ stating that
${\mathcal S}(V)$ is always Clifford \cite{BS}. In
\cite{Ba1,Ba2,Ba3}, Bazzoni examined the case of Pr\"ufer domains of
finite character, showing that these, too, have Clifford class
semigroup. In 2001, she completely resolved the problem for the
class of integrally closed domains by proving that ``{\em $R$ is an
integrally closed domain with Clifford class semigroup if and only
if $R$ is a Pr\"ufer domain of finite character}" \cite{Ba4}, Theorem
4.5. It is worth recalling that, in the series of papers
\cite{O1,O2,O3}, Olberding undertook an extensive study of (Lipman
and Sally-Vasconcelos) stability conditions which prepared the
ground to address the correlation between stability and the theory
of class semigroups.

A domain $R$ is called a PVMD (Pr\"ufer $v$-multiplication domain)
if the $v$-finite $v$-ideals form a group under the
$t$-multiplication; equivalently, if $R_{M}$ is a valuation domain
for each $t$-maximal ideal $M$ of $R$. Ideal $t$-multiplication
converts ring notions such as PID, Dedekind, B\'ezout (of finite
character), Pr\"ufer (of finite character), and integrality to UFD,
Krull, GCD  (of finite $t$-character), PVMD (of finite
$t$-character), and pseudo-integrality, respectively. Recall at this
point that the PVMDs of finite $t$-character (i.e., each proper
$t$-ideal is contained in only finitely many $t$-maximal ideals) are
exactly the Krull-type rings introduced and studied by Griffin in
1967-68 \cite{Gr1,Gr2}. Also pseudo-integrality (which should be
termed $t$-integrality) was introduced and studied in 1991 by D. F.
Anderson, Houston and Zafrullah \cite{AHZ}. We'll provide more
details about this property which turned to be crucial for our
study.

This paper examines ring-theoretic properties of an integral domain
which reciprocally reflect in semigroup-theoretic properties of its
$t$-class semigroup. Notions and contexts that suit best
$t$-multiplication are studied in an attempt to parallel analogous
developments and generalize well-known results on class semigroups.
Recall from \cite{Ba4,KM} that an integral domain $R$ is Clifford
regular (resp., Boole regular) if ${\mathcal S}(R)$ is a Clifford
(resp., Boolean) semigroup. A first correlation between regularity
and stability conditions can be sought through Lipman stability.
Indeed, $R$ is called an L-stable domain if $\bigcup_{n\geq 1}
(I^{n}:I^{n})=(I:I)$ for every nonzero ideal $I$ of $R$ \cite{AHP}.
Lipman introduced the notion of stability in the specific setting of
one-dimensional commutative semi-local Noetherian rings in order to
give a characterization of Arf rings; in this context, L-stability
coincides with Boole regularity \cite{Lip}. By analogy, we call an
integral domain $R$ Clifford (resp., Boole) $t$-regular if
${\mathcal S}_{t}(R)$ is a Clifford (resp., Boolean) semigroup.
Clearly, a Boole $t$-regular domain is Clifford $t$-regular.

Section~\ref{sec:2} establishes $t$-analogues of basic results on
$t$-regularity. We notice that a Krull domain (resp., UFD) is
Clifford (resp., Boole) $t$-regular. These two classes of domains
serve as a starting ground for $t$-regularity (as Dedekind domains
and PIDs do for regularity). We show that $t$-regularity stands as a
default measure for some classes of Krull-like domains. For
instance, it measures how far a $t$-almost Dedekind domain \cite{Kg}
is from being a Krull domain or a UFD.  In particular, we'll see
that ``UFD = Krull + Boole $t$-regular." While an integrally closed
Clifford regular domain is Pr\"ufer \cite{ZZ}, an integrally closed
Clifford $t$-regular domain need not be a PVMD; an example is
provided in this regard (Cf. Example~\ref{sec:2.7}). As a prelude to
this, our main theorem of this section (Theorem~\ref{sec:2.7})
investigates the transfer of $t$-regularity to pseudo-valuation
domains; namely, a PVD $R$ is always Clifford $t$-regular; moreover,
$R$ is Boole $t$-regular if and only if it is issued from a Boole regular valuation domain.

Section~\ref{sec:3} seeks a satisfactory $t$-analogue for Bazzoni's
theorem on Pr\"ufer domains of finite character \cite{Ba4}, Theorem
4.5 (quoted above). From \cite{AHZ}, the pseudo-integral closure of a domain $R$ is defined as $\widetilde{R}=\cup(I_{t}\colon I_{t})$, where $I$ ranges over the set of finitely generated ideals of $R$; and $R$ is said to be
pseudo-integrally closed if $R=\widetilde{R}$. Clearly $R'\subseteq
\widetilde{R}\subseteq \overline{R}$, where $R'$ and $\overline{R}$
are respectively the integral closure and the complete integral
closure of $R$. In view of Example~\ref{sec:2.7} (mentioned above),
one has to elevate the ``integrally closed" assumption in regularity
results to ``pseudo-integrally closed." In this vein, we conjecture
that ``{\em a pseudo-integrally closed domain is Clifford
$t$-regular if and only if it is a Krull-type domain}." Our main
theorem of this section (Theorem~\ref{sec:3.2}) asserts that ``{\em
a PVMD is Clifford $t$-regular if and only if it is a Krull-type
domain}." It recovers Bazzoni's theorem and also reveals the fact
that in the class of PVMDs, Clifford $t$-regularity coincides with
the finite $t$-character condition. Moreover, we are able to
validate the conjecture in a large class of integral domains
(Corollary~\ref{sec:3.12}).

Section~\ref{sec:4} is devoted to generating examples. We treat the
possible transfer of the PVMD notion endowed with the finite
$t$-character condition to pullbacks and polynomial rings. Original
families of integral domains with Clifford $t$-class
semigroup stem from our results.

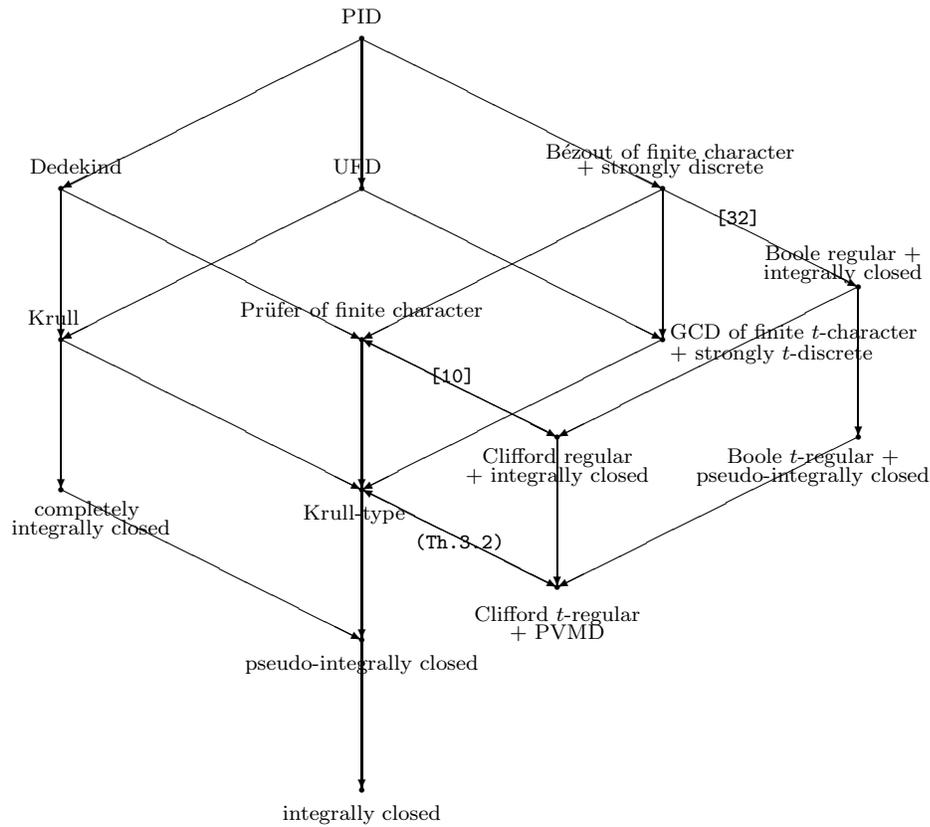
\begin{figure}[h!]
\centering
\[\setlength{\unitlength}{1mm}
\begin{picture}(100,105)(0,-134)
\put(40,-30){\vector(0,-1){20}}\put(40,-30){\vector(2,-1){40}}
\put(40,-30){\vector(-2,-1){40}}\put(0,-50){\vector(0,-1){20}}
\put(0,-50){\vector(2,-1){40}}\put(80,-50){\vector(0,-1){20}}
\put(80,-50){\vector(-2,-1){40}} 
\put(106,-63){\vector(0,-1){20}}
\put(106,-63){\vector(-2,-1){40}}
\put(66,-83){\vector(0,-1){20}}\put(106,-83){\vector(-2,-1){40}}
 \put(40,-50){\vector(2,-1){40}}
\put(40,-50){\vector(-2,-1){40}}\put(0,-70){\vector(2,-1){40}}
\put(0,-70){\vector(0,-1){20}} \put(80,-70){\vector(-2,-1){40}}

\put(40,-70){\vector(0,-1){20}}\put(40,-90){\vector(0,-1){20}}
\put(0,-90){\vector(2,-1){40}} \put(40,-110){\vector(0,-1){20}}
\put(40,-70){\vector(2,-1){26}}\put(40,-90){\vector(2,-1){26}}
\put(66,-103){\vector(-2,1){26}}\put(66,-83){\vector(-2,1){26}}
\put(80,-50){\vector(2,-1){26}}

\put(52,-76){\makebox(0,0)[b]{\footnotesize\tt\cite{Ba4}}}
\put(90,-55){\makebox(0,0)[b]{\footnotesize\tt\cite{KM}}}
\put(53,-98){\makebox(0,0)[b]{\footnotesize\tt(Th.\ref{sec:3.2})}}
\put(40,-30){\circle*{.7}} \put(40,-28){\makebox(0,0)[b]{\footnotesize
PID}} \put(0,-50){\circle*{.7}} \put(2,-48){\makebox(0,0)[b]{\footnotesize
Dedekind}} \put(40,-50){\circle*{.7}}
\put(39.4,-48){\makebox(0,0)[b]{\footnotesize UFD}}
\put(80,-50){\circle*{.7}} \put(81,-46){\makebox(0,0)[b]{\footnotesize
B\'ezout of finite character}} \put(81,-46){\makebox(0,0)[t]{\footnotesize
+ strongly discrete}}
\put(0,-70){\circle*{.7}}
\put(-1,-68){\makebox(0,0)[b]{\footnotesize Krull}}

\put(80,-70){\circle*{.7}}\put(81,-69){\makebox(0,0)[l]{\footnotesize GCD of finite $t$-character}} 

\put(81,-72){\makebox(0,0)[l]{\footnotesize + strongly $t$-discrete}}
\put(40,-70){\circle*{.7}}\put(40,-67){\makebox(0,0)[b]{\footnotesize
Pr\"ufer of finite character}}
 \put(40,-90){\circle*{.7}} \put(39,-92){\makebox(0,0)[t]{\footnotesize Krull-type}}
\put(0,-90){\circle*{.7}} \put(4,-94){\makebox(0,0)[b]{\footnotesize completely
}}\put(4,-94){\makebox(0,0)[t]{\footnotesize integrally closed}}
\put(40,-110){\circle*{.7}}
\put(40,-112){\makebox(0,0)[t]{\footnotesize pseudo-integrally closed}}
\put(40,-130){\circle*{.7}} \put(40,-132){\makebox(0,0)[t]{\footnotesize
integrally closed}} 
\put(106,-63){\circle*{.7}}\put(104,-60){\makebox(0,0)[b]{\footnotesize Boole regular +}}\put(104,-60){\makebox(0,0)[t]{\footnotesize integrally closed}}

\put(106,-83){\circle*{.7}}\put(100,-87){\makebox(0,0)[b]{\footnotesize Boole $t$-regular +}}

\put(100,-87){\makebox(0,0)[t]{\footnotesize
pseudo-integrally closed}}
\put(66,-83){\circle*{.7}}\put(66,-87){\makebox(0,0)[b]{\footnotesize
Clifford regular}}\put(66,-87){\makebox(0,0)[t]{\footnotesize
+ integrally closed}}
\put(66,-103){\circle*{.7}}\put(66,-108){\makebox(0,0)[b]{\footnotesize
Clifford $t$-regular}}\put(66,-108){\makebox(0,0)[t]{\footnotesize
+ PVMD}}
\end{picture}\]
\caption{A ring-theoretic perspective for ($t$-)regularity}\label{D1}
\end{figure}

All rings considered in this paper are integral domains. For the
convenience of the reader, Figure~\ref{D1} displays a
 diagram of implications summarizing the relations between the main
classes of integrally closed domains that provide a suitable
environment for our study. It also places ($t$-)regularity in a
ring-theoretic perspective.

\bigskip
{\bf Acknowledgments.} This work was funded by King Fahd University of Petroleum \& Minerals under Research Project \# MS/t-Class/257.
The authors would like to thank the referees of
the paper for a very careful reading and useful comments and suggestions.

\section{\bf Basic results on $t$-regularity}\label{sec:2}

Let $R$ be a domain with quotient field $K$. We first review some
terminology related to the $v$- and $t$-operations. For a nonzero
fractional ideal $I$ of $R$, let $I^{-1}$ denote $(R:I)=\{x\in K\mid
xI\subseteq R\}$. The $v$- and $t$-closures of $I$ are defined,
respectively, by $I_v=(I^{-1})^{-1}$ and $I_t=\cup J_v$ where $J$
ranges over the set of finitely generated subideals of $I$. The
(nonzero) ideal $I$ is said to be divisorial or a $v$-ideal if
$I_v=I$, and a $t$-ideal if $I_t=I$. Under the ideal
$t$-multiplication $(I,J)\mapsto (IJ)_t$, the set $F_{t}(R)$ of
fractional $t$-ideals of $R$ is a semigroup with unit $R$. An
invertible element for this operation is called a $t$-invertible
$t$-ideal of $R$. So that the set $Inv_{t}(R)$ of $t$-invertible
fractional $t$-ideals of $R$ is a group with unit $R$. For more
basic details about star operations, we refer the reader to
\cite{Gi}, Sections 32 and 34. Let $F(R)$, $Inv(R)$, and $P(R)$
denote the sets of nonzero, invertible, and nonzero principal
fractional ideals of $R$, respectively. Under this notation, the
($t$-)class groups and semigroups are defined as follows:
$\Pic(R)=Inv(R)/P(R)$, $\Cl(R)=Inv_{t}(R)/P(R)$, ${\mathcal
S}(R)=F(R)/P(R)$, and ${\mathcal S}_{t}(R)=F_{t}(R)/P(R)$.

Recall two basic properties of the $t$-operation which will be used
(in different forms) throughout the paper. For any two nonzero
ideals $I$ and $J$ of a domain $R$, we have
$(IJ)_{t}=(I_{t}J)_{t}=(IJ_{t})_{t}=(I_{t}J_{t})_{t}$. Also one can
easily check that $(I_{t}:J)=(I_{t}:J_{t})$. In particular, we have
$I^{-1}=(R:I)=(R:I_{t})$ and, if $I$ is a $t$-ideal,
$(I:I^{2})=(I:(I^{2})_{t})$. Actually, these properties hold for any
star operation.

Throughout, we shall use $\qf(R)$ to denote the quotient field of a
domain $R$ and ${\overline I}$ to denote the isomorphy class of an
ideal $I$ of $R$ in ${\mathcal S}_{t}(R)$.\bigskip

Our first result displays necessary and/or sufficient
ideal-theoretic conditions for the isomorphy class of an ideal to be
regular in the $t$-class semigroup.

\begin{lemma}\label{sec:2.1}
Let $I$ be a $t$-ideal of a domain $R$. Then:\\
\1 ${\overline I}$ is regular in ${\mathcal S_t}(R)$ if and only if
$I=(I^{2}(I:I^{2}))_{t}$.\\
\2 If $I$ is $t$-invertible, then ${\overline I}$ is regular in
${\mathcal S_t}(R)$.
\end{lemma}

\begin{proof}
\1  Assume ${\overline I}$ is regular in ${\mathcal S_t}(R)$. Then
there exist a  fractional $t$-ideal $J$ of $R$ and $0\not =c\in
\qf(R)$ such that $I = c(JI^{2})_{t} = (cJI^{2})_{t}$. We may denote
$cJ$ by $J$, that is, $I = (JI^{2})_{t}$. Since $JI^{2}\subseteq
(JI^{2})_{t} = I$, $J\subseteq (I:I^{2})$. So
$I=(JI^{2})_{t}\subseteq (I^{2}(I:I^{2}))_{t}\subseteq I$ and hence
$I=(I^{2}(I:I^{2}))_{t}$. The converse is trivial.

\2 Assume $(II^{-1})_{t} = R$. Then $R\subseteq (I:I)\subseteq
(II^{-1}:II^{-1})\subseteq ((II^{-1})_{t}:(II^{-1})_{t}) =R$.  So
$(I:I^{2}) = ((I:I):I)=I^{-1}$. Hence $(I^{2}(I:I^{2}))_{t} =
(I^{2}I^{-1})_{t}= (I(II^{-1}))_{t}=(I(II^{-1})_{t})_{t}=I$. By \1,
${\overline I}$ is regular in ${\mathcal S_t}(R)$. \end{proof}

Next, we show that Krull domains (resp., UFDs) are Clifford (resp.,
Boole) $t$-regular. Further, we identify $t$-regularity as a default
condition for some classes of Krull-like domains towards the Krull
(or UFD) property. Recall at this point that a domain $R$ is Krull if and only if every $t$-ideal of $R$ is $t$-invertible.

\begin{proposition}\label{sec:2.2}
\1  Any Krull domain is Clifford $t$-regular.\\
\2 A domain $R$ is a UFD if and only if $R$ is Krull and Boole
$t$-regular.
\end{proposition}

\begin{proof}
\1  Follows from Lemma~\ref{sec:2.1}(2).

\2 Clearly, a UFD is Boole $t$-regular. We need only prove the ``if"
assertion. Assume $R$ is Krull and Boole $t$-regular and let $I$ be
a $t$-ideal of $R$. There exists $0\not=c\in \qf(R)$ such that
$(I^{2})_{t}=cI$. Then
$(I:I^{2})=(I:(I^{2})_{t})=(I:cI)=c^{-1}(I:I)$. $R$ is completely
integrally closed, then $(I:I)=R$, so that $(I:I^{2})=(R:I)=I^{-1}$.
Therefore $I^{-1}=c^{-1}R$, and hence $II^{-1}=c^{-1}I$. Since $I$
is $t$-invertible, $R=(II^{-1})_{t}=(c^{-1}I)_{t}=c^{-1}I$, hence
$I=cR$. It follows that $\Cl(R)={\mathcal S_t}(R)=0$, i.e., $R$ is a
UFD. \end{proof}

Recall from \cite{Kg} that a domain $R$ is said to be $t$-almost
Dedekind if $R_{M}$ is a rank-one DVR for each $t$-maximal ideal $M$
of $R$. This notion falls strictly between the classes of Krull
domains and PVMDs. Our next result shows that $t$-regularity
measures how far a $t$-almost Dedekind domain or completely
integrally closed domain is from being Krull or a UFD. A domain $R$
is said to be strongly $t$-discrete if it has no $t$-idempotent
$t$-prime ideals, i.e., for every  $t$-prime ideal $P$ of $R$,
$(P^{2})_{t}\subsetneq P$ (Cf. \cite{E}).

\begin{proposition}\label{sec:2.3}
Let $R$ be a domain. The following statements are equivalent:\\
\i $R$ is Krull (resp., a UFD);\\
\ii $R$ is $t$-almost Dedekind and Clifford (resp., Boole) $t$-regular;\\
\iii $R$ is strongly $t$-discrete, completely integrally closed, and
Clifford (resp., Boole) $t$-regular.
\end{proposition}

\begin{proof}
\i $\Longrightarrow$ \ii Straightforward.

\ii $\Longrightarrow$ \i Suppose there exists a  $t$-ideal $I$ of
$R$ which is not $t$-invertible. Then $J=(II^{-1})_{t}$ is a proper
trace $t$-ideal of $R$ with $J^{-1} = (J:J)$. Further, $R$ is
completely integrally closed since $R = \bigcap R_{M}$, where $M$
ranges over the $t$-maximal ideals of $R$ \cite{Kg}, Proposition
2.9. Therefore $J^{-1}=(J:J)=R$, so that $J^{2}(J:J^{2}) =
J^{2}((J:J):J)=J^{2}J^{-1}= J^{2}$. Now, ${\overline J}$ is regular
in ${\mathcal S_t}(R)$, then $ J=(J^{2}(J:J^{2}))_{t} =
(J^{2})_{t}$. By induction, we get $J= (J^{n})_{t}$, for each $n\geq
1$. By \cite{Kg}, Proposition 2.54, $J=\bigcap_{n\geq 1}
(J^{n})_{t}=(0)$, the desired contradiction.

\i $\Longrightarrow$ \iii Let $P$ be a $t$-prime ideal of $R$. Since
$R$ is Krull, $(PP^{-1})_{t}=R$. Suppose $P$ is $t$-idempotent,
i.e., $(P^{2})_{t}=P$. Then
$((P^{2})_{t}P^{-1})_{t}=(PP^{-1})_{t}=R$. Hence
$P=((PP^{-1})_{t}P)_{t}=
(P^{2}P^{-1})_{t}=((P^{2})_{t}P^{-1})_{t}=R$, absurd.

\iii $\Longrightarrow$ \i Suppose there is a  $t$-ideal $I$ of $R$
such that $J=(II^{-1})_{t}\subsetneq R$. Here too we have
$J^{-1}=(J:J)=R$. Let $M$ be a $t$-maximal ideal of $R$ containing
$J$. Necessarily, $(M:M)=M^{-1}=R$. Therefore
$(M:M^{2})=((M:M):M)=(R:M)=M^{-1}=R$. So $M^{2}(M:M^{2})=M^{2}$.
Since $R$ is Clifford $t$-regular, then
$M=(M^{2}(M:M^{2}))_{t}=(M^{2})_{t}$ and hence $M$ is
$t$-idempotent, absurd.

 The Boolean statements follow readily from the
Clifford statements combined with Proposition~\ref{sec:2.2},
completing the proof. \end{proof}

Notice that the ring of all entire functions in the complex plane is
(B\'ezout) strongly ($t$-)discrete  \cite{FHP}, Corollary 8.1.6, and
completely integrally closed, but it is not ($t$-)almost Dedekind
(since it has an infinite Krull dimension). Also the ``strongly
$t$-discrete" assumption in \iii is not superfluous, since a
non-discrete rank-one valuation domain is completely integrally
closed and Clifford ($t$-)regular \cite{BS}, but it is not Krull.
\bigskip

The next result establishes the transfer of $t$-regularity to
polynomial rings. Recall at this point that Clifford or Boole
regularity of a polynomial ring $R[X]$ forces $R$ to be a field
\cite{KM}, Corollary 2.5.

\begin{proposition}\label{sec:2.4}
Let $R$ be an integrally closed domain and $X$ an indeterminate over
$R$. Then $R$ is Clifford (resp., Boole) $t$-regular if and only if
so is $R[X]$.\end{proposition}

\begin{proof} Assume that $R$ is Clifford $t$-regular and let $J$ be a $t$-ideal of
$R[X]$ with $I= J\cap R$. If $I\not= 0$, then $I$ is a $t$-ideal of
$R$ and hence $J = I[X]$. If $I=(0)$, then $J=fA[X]$ for some $f\in
R[X]$ and $A$ a fractional $t$-ideal  of $R$ \cite{Q}. So that
$J^{2}(J:J^{2})$ equals $(I^{2}(I:I^{2}))[X]$ or
$f(A^{2}(A:A^{2}))[X]$. In both cases, $(J^{2}(J:J^{2}))_{t} =J$ by
\cite{Kg}, Proposition 2.3(1) (which ensures that the $t$-operation
is stable under ideal extension). Therefore ${\overline J}$ is
regular in ${\mathcal S_t}(R[X])$. Conversely, If $I$ is a $t$-ideal
of $R$, consider the $t$-ideal $I[X]$ of $R[X]$ and apply the same
techniques backward. Similar arguments as above lead to the
conclusion for the Boolean statement. \end{proof}

The next result establishes the transfer of $t$-regularity to two
types of overrings.

\begin{proposition}\label{sec:2.5}
Let $R$ be a Clifford (resp., Boole) $t$-regular domain. Then:\\
\1 $R_{S}$ is Clifford (resp., Boole) $t$-regular, for any multiplicative subset $S$ of $R$.\\
\2 $(I_{v}:I_{v})$ is Clifford (resp., Boole) $t$-regular, for any
nonzero ideal $I$ of $R$.
\end{proposition}

For the proof, we need the following lemma.

\begin{lemma}\label{sec:2.5.1}
Let $R$ be a domain, I a fractional ideal of $R$, and $S$ a
multiplicative subset of $R$. Then $I_{t}\subseteq
(IR_{S})_{t_{1}}$, where $t_{1}$ denotes the $t$-operation with
respect to $R_{S}$.
\end{lemma}

\begin{proof}
 Let $x\in I_{t}$. Then there exists a finitely generated ideal $A$ of
$R$ such that $A\subseteq I$ and $x(R:A)\subseteq R$. Hence
$x(R_{S}:AR_{S})=x(R:A)R_{S}\subseteq R_{S}$. Therefore $x\in
(AR_{S})_{t_{1}}\subseteq (IR_{S})_{t_{1}}$. \end{proof}

\begin{proof}[Proof of Proposition~\ref{sec:2.5}]
\1  If $J$ is a $t$-ideal of $R_{S}$, then $I=J\cap R$ is a
$t$-ideal of $R$ by Lemma~\ref{sec:2.5.1}. Since $R$ is Clifford
(resp., Boole) $t$-regular, then $I=(I^{2}(I:I^{2}))_{t}$ (resp.,
$(I^{2})_{t}=cI$ for some nonzero $c\in \qf(R)$). Hence
$J=IR_{S}=(I^{2}(I:I^{2}))_{t}R_{S}\subseteq
((I^{2}(I:I^{2}))R_{S})_{t_{1}}\subseteq
(J^{2}(J:J^{2}))_{t_{1}}\subseteq J$ (resp.,
$cJ=cIR_{S}=(I^{2})_{t}R_{S}\subseteq
(I^{2}R_{S})_{t_{1}}=(J^{2})_{t_{1}}\subseteq cJ$, since
$I^{2}\subseteq (I^{2})_{t}=cI$ and then $J^{2}\subseteq cJ$).
Therefore $J=(J^{2}(J:J^{2}))_{t_{1}}$ (resp.,
$(J^{2})_{t_{1}}=cJ$). It follows that $R_{S}$ is Clifford (resp.,
Boole) $t$-regular.

\2 Let $I$ be a nonzero ideal of $R$ and set $T=(I_{v}:I_{v})$.
Since $T=(II^{-1})^{-1}=(II^{-1}:II^{-1})=
((II^{-1})_{v}:(II^{-1})_{v})$, without loss of generality, we may
assume that $I$ is a trace $v$-ideal of $R$, that is
$T=I^{-1}=(I:I)$. Also denote by $v_{1}$ and $t_{1}$ the $v$- and
$t$- operations with respect to $T$. Let $J$ be a nonzero ideal of
$T$. Then $J$ is a fractional ideal of $R$ and we claim that
$J_{t}\subseteq J_{t_{1}}$. Indeed, let $x\in J_{t}$. Then there
exists a finitely generated (fractional) ideal $A$ of $R$ such that
$A\subseteq J$ and $x(R:A)\subseteq R$. Let $z\in (T:AT)$. Then
$zAI\subseteq I\subseteq R$, hence $zI\subseteq (R:A)$, whence
$xzI\subseteq x(R:A)\subseteq R$ and $xz\in I^{-1}=T$. Therefore
$x(T:AT)\subseteq T$, and hence $x\in
(AT)_{v_{1}}=(AT)_{t_{1}}\subseteq J_{t_{1}}$. Consequently, if $J$
is a $t$-ideal of $T$, then it's a $t$-ideal of $R$. Since $R$ is
Clifford (resp., Boole) $t$-regular, then
$J=(J^{2}(J:J^{2}))_{t}\subseteq (J^{2}(J:J^{2}))_{t_{1}}\subseteq
J$ (resp., $cJ=(J^{2})_{t}\subseteq (J^{2})_{t_{1}}\subseteq cJ$,
since $J^{2}\subseteq (J^{2})_{t}=cJ$, for some nonzero $c\in
\qf(R)=\qf(T)$). Hence $J=(J^{2}(J:J^{2}))_{t_{1}}$ (resp.,
$(J^{2})_{t_{1}}=cJ$) and therefore $T$ is Clifford (resp., Boole)
$t$-regular. \end{proof}

We close this section with an investigation of the integrally closed
setting. In this vein, recall Zanardo-Zannier's crucial result
 that an integrally closed Clifford regular domain is necessarily
Pr\"ufer \cite{ZZ}. In \cite{KM}, we stated an analogue for Boole
regularity, that is, an integrally closed Boole regular domain is
B\'ezout. Next, we show that an integrally closed Clifford (or Boole)
$t$-regular domain need not be a PVMD, the natural context for
$t$-regularity. Our family of such examples stems from the following
theorem on the inheritance of $t$-regularity by PVDs (i.e.,
pseudo-valuation domains). We refer the reader to \cite{HH} for the definition and the main properties of PVDs.

\begin{thm}\label{sec:2.6}
Let $R$ be a PVD. Then:\\
\1  $R$ is Clifford $t$-regular.\\
\2 $R$ is Boole $t$-regular if and only if its associated valuation
overring is Boole regular.
\end{thm}

\begin{proof} \1  We may assume that $R$ is not a valuation domain.
Proposition 2.6 of \cite{AD} characterizes PVDs in terms of
pullbacks. The aforementioned proposition states that $R$ is a PVD
if and only if $R=\phi^{-1}(k)$ for some subfield $k$ of $K=V/M$,
where $V$ is the associated valuation overring of $R$, $M$ its
maximal ideal and $\phi$ the canonical homomorphism from $V$ onto
$K$. Now, let $I$ be a $t$-ideal of $R$. If $I$ is an ideal of $V$,
we are done (since $V$ is Clifford regular). If $I$ is not an ideal
of $V$, then $I=c\phi^{-1}(W)$, where $0\not =c\in M$ and $W$ is a
$k$-vector space such that $k\subseteq W\subset K$ (Cf.
\cite{BG}, Theorem 2.1(n)). Assume $k\subsetneq W$. Then
$(k:W)=(0)$. Hence $I^{-1}=(R:I) =
(\phi^{-1}(k):c\phi^{-1}(W))=c^{-1}\phi^{-1}(k:W)=c^{-1}M$ by \cite{HKLM1},
Proposition 6. Since $R$ is a PVD which is not a valuation
domain, by \cite{HZ}, Proposition 4.3, $R$ is a $TV$-domain (i.e.
the $t$- and $v$- operations coincide in $R$). Hence
$I=I_{t}=I_{v}=(R:c^{-1}M)=cM^{-1}= cV$ is an ideal of $V$, a
contradiction. Therefore $k=W$ and then $I=cR$ is a principal ideal
of $R$. So ${\overline I}$
is regular in ${\mathcal S_t}(R)$, as desired.\\
\2 Assume that $R$ is Boole $t$-regular. By
Proposition~\ref{sec:2.5}, $V=(M:M)=(M_{v}:M_{v})$ is Boole regular
(the $t$-operation on $V$ is trivial). Conversely, assume that $V$ is Boole regular.
Similar arguments as above lead to the conclusion.\end{proof}

Contrast this result with \cite{KM}, Theorem 5.1, which asserts that
a PVD $R$ associated to a valuation (resp., strongly discrete
valuation) domain $(V, M)$ is Clifford (resp., Boole) regular if and
only if $[V/M:R/M]=2$.

\begin{example}\label{sec:2.7}\rm
There exists an integrally closed Boole (hence Clifford) $t$-regular
domain which is not a PVMD. Indeed, let $k$ be a field and let $X$
and $Y$ be two indeterminates over $k$. Let $R=k+M$ be the PVD
associated to the rank-one DVR $V=k(X)[[Y]]=k(X)+M$, where $M=YV$.
Clearly, $R$ is integrally closed and, by Theorem~\ref{sec:2.6}, $R$
is Boole $t$-regular. However, $R$ is not a PVMD by \cite{FG}, Theorem
4.1.
\end{example}

\section{\bf Clifford $t$-regularity}\label{sec:3}

Recall from \cite{Gr2} that a Krull-type domain is a PVMD with
finite $t$-character (i.e., each nonzero nonunit is contained in
only finitely many $t$-maximal ideals). Also a domain $R$ is said to
be pseudo-integrally closed if $R=\widetilde{R}=\cup(I_{t}\colon
I_{t})$, where $I$ ranges over the set of finitely generated ideals
of $R$ \cite{AHZ}. This section seeks a $t$-analogue for Bazzoni's
theorem that ``an integrally closed domain $R$ is Clifford regular
if and only if $R$ is a Pr\"ufer domain of finite character"
\cite{Ba4}, Theorem 4.5. In view of Example~\ref{sec:2.7}, one has
to elevate the ``integrally closed" assumption to
``pseudo-integrally closed." Accordingly, we claim the following:

\begin{conjecture}\label{sec:3.1}
A pseudo-integrally closed domain $R$ is Clifford $t$-regular if and
only if $R$ is a Krull-type domain.
\end{conjecture}

This is still elusively open. Yet, our main result (Theorem
~\ref{sec:3.2}) of this section recovers Bazzoni's theorem and
validates this conjecture in large classes of (pseudo-integrally
closed) domains.

\begin{thm}\label{sec:3.2}
A PVMD is Clifford $t$-regular if and only if it is a Krull-type
domain {\rm (}i.e., in a PVMD, Clifford $t$-regularity coincides
with the finite $t$-character condition{\rm )}.
\end{thm}

The proof of the theorem involves several preliminary results, some
of which are of independent interest. Experts of $t$- operation may
skip the proofs of Lemmas \ref{sec:3.8}, \ref{sec:3.9} and
\ref{sec:3.10} which are similar in form to their respective
analogues for the trivial operation.

The following notation, connected with the $t$-ideal structure of a
PVMD, will be of use in the sequel. Assume $R$ is a PVMD and let $I$
be a  $t$-ideal of $R$ and $x$ a nonzero nonunit element of $R$. We
shall use $\Max_{t}(R)$ to denote the set of maximal $t$-ideals of
$R$. Set $\Max_{t}(R, I)=\{M\in \Max_{t}(R)\ | \ I\subseteq M\}$,
$\Max_{t}(R, x)=\Max_{t}(R, xR)$, and ${\mathcal{T}}_{t}(R)=\{M\in
\Max_{t}(R)\ | \ R_{M}\not\supseteq\bigcap_{M\not=N}R_{N}, N\in
\Max_{t}(R)\}$. Finally, given $M$ and $N$ two $t$-maximal ideals of
$R$, we denote by $M\wedge N$ the largest prime ideal of $R$
contained in $M\cap N$. Note that prime ideals of $R$ contained in a
$t$-maximal ideal are necessarily $t$-ideals and form a chain.

\begin{lemma}\label{sec:3.3}
Let $R$ be a PVMD and $I$ a fractional ideal of $R$. Then for every
$t$-prime ideal $P$ of $R$, $I_{t}R_{P}=IR_{P}$.
\end{lemma}

\begin{proof} Here $R_{P}$ is a valuation domain (where the $t$- and trivial operations coincide), so Lemma~\ref{sec:2.5.1}
leads to the conclusion. \end{proof}

\begin{lemma}\label{sec:3.4}
Let $R$ be a PVMD which is Clifford $t$-regular and $I$ a
 nonzero fractional ideal of $R$.  Then $I$ is $t$-invertible
if and only if $I$ is $t$-locally principal.
\end{lemma}

\begin{proof}  Suppose $I$ is $t$-locally principal and set $J=(II^{-1})_{t}$.
Let $M\in \Max_{t}(R)$. Then $IR_{M}=aR_{M}$ for some nonzero $a\in
I$. By Lemma~\ref{sec:3.3}, $(I_{t}:I_{t})\subseteq
(I_{t}R_{M}:I_{t}R_{M})=(IR_{M}:IR_{M})=(aR_{M}:aR_{M})=R_{M}$.
Therefore $R\subseteq(I_{t}:I_{t})\subseteq\bigcap_{M\in
\Max_{t}(R)}R_{M}=R$. So
$(I_{t}:I_{t}^{2})=((I_{t}:I_{t}):I_{t})=(R:I_{t})=I^{-1}$. Since
$R$ is Clifford $t$-regular, then
$I_{t}=(I_{t}^{2}(I_{t}:I_{t}^{2}))_{t}=(I_{t}^{2}I^{-1})_{t}=(IJ)_{t}$.
By Lemma~\ref{sec:3.3},
$aR_{M}=IR_{M}=I_{t}R_{M}=(IJ)_{t}R_{M}=IJR_{M}=aJR_{M}$. It follows
that $R_{M}=JR_{M}$ for every $M\in \Max_{t}(R)$, which forces $J$
to equal $R$, as desired.

Conversely, assume that $I$ is $t$-invertible. Then there is a
finitely generated ideal $J$ of $R$ such that $J\subseteq I$ and
$I_{t}=J_{t}$. Hence for each $M\in \Max_{t}(R)$,
$IR_{M}=I_{t}R_{M}=J_{t}R_{M}=JR_{M}=aR_{M}$ for some $a\in J$,
since $R_{M}$ is a valuation domain. \end{proof}

\begin{lemma}\label{sec:3.5}
Let $R$ be a PVMD which is Clifford $t$-regular and let $P\subsetneq
Q$ be two $t$-prime ideals of $R$. Then there exists a finitely
generated ideal $I$ of $R$ such that $P\subsetneq I_{t}\subseteq Q$.
\end{lemma}

\begin{proof} Let $x\in Q\setminus P$ and set $J= xR+P$. We claim that $J$ is $t$-invertible.
Indeed, let $M$ be a $t$-maximal ideal of $R$. If $J\not\subseteq
M$, then $JR_{M}=R_{M}$. If $J\subseteq M$, then $PR_{M}\subsetneq
xR_{M}$ since $R_{M}$ is a valuation domain, whence $JR_{M}=xR_{M}$.
So $J$ is $t$-locally principal and hence $t$-invertible by
Lemma~\ref{sec:3.4}. Therefore $J_{t}=I_{t}$ for some finitely
generated subideal $I$ of $J$. It follows that $P\subsetneq
I_{t}\subseteq Q$. \end{proof}

\begin{lemma}\label{sec:3.6}
Let $R$ be a PVMD which is Clifford $t$-regular and $P$ a $t$-prime
ideal of $R$. Then $E(P)=(P:P)$ is a PVMD which is Clifford
$t$-regular and $P$ is a $t$-maximal ideal of $E(P)$.
 \end{lemma}

\begin{proof} If $P\in \Max_{t}(R)$, then $E(P)=R$.
We may then assume that $P\notin \Max_{t}(R)$. By \cite{H}, Proposition
1.2 and Lemma 2.4, $E(P)=P^{-1}$ and $P$ is a $t$-prime ideal of
$E(P)$. Further $E(P)=P^{-1}$ is $t$-linked over $R$, so $E(P)$ is a
PVMD \cite{Kg}. Let $t_{1}$ and $v_{1}$ denote the $t$- and $v$-
operations with respect to $E(P)$ and let $J$ be a nonzero
fractional ideal of $E(P)$. Clearly $J$ is a fractional ideal of $R$
and we claim that $J_{t}\subseteq J_{t_{1}}$. Indeed, let $x\in
J_{t}$. Then there is a finitely generated subideal $I$ of $J$ such
that $x\in I_{v}$. So $xI^{-1}\subseteq R$. Let $z\in
(E(P):IE(P))=(P^{-1}:IP^{-1})$. Then $zIP^{-1}\subseteq P^{-1}$. So
$zIP\subseteq P\subseteq R$. Then $zP\subseteq I^{-1}$. So
$xzP\subseteq xI^{-1}\subseteq R$. Hence $xz\in P^{-1}=E(P)$. So
$x(E(P):IE(P))\subseteq E(P)$ and therefore $x\in
(IE(P))_{v_{1}}\subseteq J_{t_{1}}$. Now, let $J$ be a $t$-ideal of
$E(P)$. By the above claim $J$ is a $t$-ideal of $R$. Since $R$ is
Clifford $t$-regular, then $J=(J^{2}(J:J^{2}))_{t}\subseteq
(J^{2}(J:J^{2}))_{t_{1}}\subseteq J$ and therefore
$J=(J^{2}(J:J^{2}))_{t_{1}}$. It follows that $E(P)$ is Clifford
$t$-regular. To complete the proof, we need to show that $P$ is a
$t$-maximal ideal of $E(P)$. Deny. Then there is a $t$-maximal ideal
$Q$ of $E(P)$ such that $P\subsetneq Q$. By Lemma~\ref{sec:3.5},
there is a finitely generated ideal $J$ of $E(P)$ such that
$P\subsetneq J_{t_{1}}\subseteq Q$. On the other hand, since $P$ is
a non-$t$-maximal $t$-prime ideal of $E(P)$, by \cite{H}, Proposition
1.2, $(E(P):P)=(P:P)=E(P)$. It follows that
$E(P)=P_{v_{1}}\subseteq
(J_{t_{1}})_{v_{1}}=J_{v_{1}}=J_{t_{1}}\subseteq Q$, the desired
contradiction. \end{proof}

\begin{lemma}\label{sec:3.7}
Let $R$ be a PVMD which is Clifford $t$-regular and $Q$ a $t$-prime
ideal of $R$. Suppose there is a nonzero prime ideal $P$ of $R$ such
that $P\subsetneq Q$ and $\htt(Q/P)=1$. Then there exists a finitely
generated subideal $I$ of $Q$ such that $\Max_{t}(R, I)=\Max_{t}(R,
Q)$.
 \end{lemma}

\begin{proof}
By \cite{Kg}, Corollary 2.47, $P$ is a $t$-prime ideal of $R$. By
\cite{H}, Proposition 1.2 and Lemma 2.4, $E(P)=P^{-1}$ and $P$ is a
$t$-prime ideal of $E(P)$. Therefore $E(Q)=(Q:Q)\subseteq
Q^{-1}\subseteq P^{-1}=E(P)$, and hence $P$ is a prime ideal of
$E(Q)$. By Lemma~\ref{sec:3.6}, $E(Q)$ is a PVMD which is Clifford
$t$-regular and $Q$ is a $t$-maximal ideal of $E(Q)$. Thus $P$ is a
$t$-prime ideal of $E(Q)$. By Lemma~\ref{sec:3.5}, there is a
finitely generated subideal $J=\sum_{1\leq i\leq n}a_{i}E(Q)$ of $Q$
such that $P\subsetneq J_{t_{1}}\subseteq Q$. We claim that
$\Max_{t}(E(Q), J)=\{Q\}$. Indeed, if there is a $t$-maximal ideal
$N$ of $E(Q)$ such that $J\subseteq N$, then $P\subsetneq
J_{t_{1}}\subseteq N$. So $P\subsetneq J_{t_{1}}\subseteq Q\wedge
N\subsetneq Q$, a contradiction since $\htt(Q/P)=1$. Now set
$I=\sum_{1\leq i\leq n}a_{i}R$. Clearly $I\subseteq Q$ and
$IE(Q)=J$. We claim that $\Max_{t}(R, I)=\Max_{t}(R, Q)$. Let $N\in
\Max_{t}(R, I)$. If $Q\not\subseteq N$, then $R\subseteq
E(Q)\subseteq R_{N}$. Hence $E(Q)_{R\setminus N}=R_{N}$. So $R_{N}$
is $t$-linked over $E(Q)$. Since $R_{N}$ is a valuation domain, then
$NR_{N}$ is a $t$-prime ideal of $R_{N}$. Hence $M=NR_{N}\cap E(Q)$
is a $t$-prime ideal of $E(Q)$. Since $I\subseteq N$, then
$I\subseteq M$. Hence $J=IE(Q)\subseteq M$. Necessarily $M\subseteq
Q$ since $\Max_{t}(E(Q), J)=\{Q\}$. So $N=NR_{N}\cap R = NR_{N}\cap
E(Q)\cap R=M\cap R\subseteq Q\cap R=Q$, which is absurd. Hence
$Q\subseteq N$ and therefore $\Max_{t}(R, I)\subseteq\Max_{t}(R,
Q)$. The reverse inclusion is trivial. \end{proof}

\begin{lemma}\label{sec:3.8}
Let $R$ be a PVMD which is Clifford $t$-regular and $M$ a
$t$-maximal ideal of $R$. If $M\in {\mathcal{T}}_{t}(R)$, then there
exists a finitely generated subideal $I$ of $M$ with $\Max_{t}(R,
I)=\{M\}$.
\end{lemma}

\begin{proof} Assume that $M\in {\mathcal{T}}_{t}(R)$. Let $x\in (\bigcap_{M\not=N\in \Max_{t}(R)} R_{N})\setminus R_{M}$. Since
$R_{M}$ is a valuation domain, then $x^{-1}\in MR_{M}$. Set
$I=x^{-1}R\cap R$. We claim that $I$ is $t$-invertible. By
Lemma~\ref{sec:3.4}, it suffices to check that $I$ is $t$-locally
principal. Let $Q$ be a $t$-maximal ideal of $R$. If $Q\not =M$,
then $I\not\subseteq Q$. Indeed, since $Q\not =M$, then $x\in
R_{Q}$. Hence $IR_{Q}=(x^{-1}R\cap R)R_{Q}=x^{-1}R_{Q}\cap
R_{Q}=R_{Q}$. So $I\not\subseteq Q$. Then $M$ is the unique
$t$-maximal ideal of $R$ that contains $I$ and $IR_{M}=x^{-1}R_{M}$,
as desired. Hence $I$ is $t$-invertible. So there is a finitely
generated subideal $J$ of $I$ such that $I_{t}=J_{t}$ and clearly
$\Max_{t}(R, J)=\{M\}$. \end{proof}

\begin{lemma}\label{sec:3.9}
Let $R$ be a PVMD which is Clifford $t$-regular. Then every nonzero
nonunit element of $R$ belongs to a finite number of $t$-maximal
ideals in ${\mathcal{T}}_{t}(R)$.
\end{lemma}

\begin{proof} Let $x$ be a nonzero nonunit element of $R$ and let $\{M_{\alpha}\}_{\alpha\in \Omega}$
be the set of all $t$-maximal ideals in ${\mathcal{T}}_{t}(R)$ that
contain $x$. For each $\alpha\in \Omega$, let $A_{\alpha}$ be a
finitely generated subideal of $M_{\alpha}$ such that $\Max_{t}(R,
A_{\alpha})=\{M_{\alpha}\}$ (Lemma~\ref{sec:3.8}). Without loss of
generality, we may assume that $x\in A_{\alpha}$ (otherwise, we
consider $B_{\alpha}=xR+A_{\alpha}$). Let $B=\sum_{\alpha\in
\Omega}(A_{\alpha})^{-1}$. Clearly, $x(A_{\alpha})^{-1}\subseteq R$,
for each $\alpha\in \Omega$. Then $xB\subseteq R$ so that $B$ is a
fractional ideal of $R$. We claim that $B$ is $t$-locally principal.
Indeed, let $N$ be a $t$-maximal ideal of $R$. We envisage two
cases. {\bf Case 1:} $N\not =M_{\alpha}$ for each $\alpha\in
\Omega$. Since $A_{\alpha}$ is finitely generated,
$(A_{\alpha})^{-1}R_{N}=(A_{\alpha}R_{N})^{-1}=R_{N}$. So
$BR_{N}=R_{N}$. {\bf Case 2:} $N=M_{\alpha}$ for some $\alpha\in
\Omega$. Then $(A_{\beta})^{-1}R_{N}=(A_{\beta}R_{N})^{-1}=R_{N}$,
for each $\beta\not = \alpha$ in $\Omega$. Hence
$BR_{N}=(A_{\alpha})^{-1}R_{N}=(A_{\alpha}R_{N})^{-1}=a^{-1}R_{N}$
where $A_{\alpha}R_{N}=aR_{N}$ (since $R_{N}$ is a valuation
domain). It follows that $B$ is $t$-invertible (Lemma~\ref{sec:3.4})
and hence there is a finitely generated subideal $J$ of $B$ such
that $J_{v}=J_{t}= B_{t}=B_{v}$. So $B^{-1}=J^{-1}$. Since $J$ is
finitely generated, then there are ${\alpha}_{1}, \dots,
{\alpha}_{r}$ such that $J\subseteq \sum_{1\leq i\leq
r}(A_{\alpha_{i}})^{-1}\subseteq B$. Therefore $B^{-1}=J^{-1}=
\big(\sum_{1\leq i\leq
r}(A_{\alpha_{i}})^{-1}\big)^{-1}=\bigcap_{1\leq i\leq
r}(A_{\alpha_{i}})_{v} =\bigcap_{1\leq i\leq
r}(A_{\alpha_{i}})_{t}$. Consequently, for each $\alpha\in \Omega$,
we have $\bigcap_{1\leq i\leq r}(A_{\alpha_{i}})_{t}=B^{-1}\subseteq
(A_{\alpha})_{v}=(A_{\alpha})_{t}\subseteq M_{\alpha}$. So there is
$\alpha_{i}$ such that $(A_{\alpha_{i}})_{t}\subseteq M_{\alpha}$,
hence $M_{\alpha}= M_{\alpha_{i}}$, whence $\alpha=\alpha_{i}$.
Therefore $\Omega=\{{\alpha}_{1}, \dots, {\alpha}_{r}\}$, as
desired. \end{proof}

\begin{lemma}\label{sec:3.10}
Let $R$ be a PVMD which is Clifford $t$-regular and $M$ a
$t$-maximal ideal of $R$. Then $M\in {\mathcal{T}}_{t}(R)$ if and
only if $M\supsetneqq\bigcup_{N}M\wedge N$ where $N$ ranges over
$\Max_{t}(R)\setminus \{M\}$
\end{lemma}

\begin{proof} Let $M\in {\mathcal{T}}_{t}(R)$ and let $A=\sum_{1\leq i\leq r}a_{i}R$
be a finitely generated subideal of $M$ such that $\Max_{t}(R,
A)=\{M\}$ (Lemma~\ref{sec:3.8}). Suppose that $M=\bigcup_{N}M\wedge
N$, where $N$ ranges over $\Max_{t}(R)\setminus \{M\}$. Then for
each $a_{i}\in A$, there is a $t$-maximal ideal $N_{i}\not =M$ such
that $a_{i}\in M\wedge N_{i}$. Since $\{M\wedge N_{i} |\ i=1, \dots,
r\}$ is a chain, let $M\wedge N_{j}$ be the largest one for some
$j\in \{1, \dots, r\}$. So $A\subseteq N_{j}$ and then $N_{j}\in
\Max_{t}(R, A)=\{M\}$, absurd.

Conversely, let $x\in M\setminus \bigcup_{N}M\wedge N$. Then, for
each $t$-maximal ideal $N\not =M$, $x^{-1}\in R_{N}$ (since $R_{N}$
is a valuation domain), hence $x^{-1}\in \bigcap_{M\not=N\in
\Max_{t}(R)} R_{N}$. Since
 $x^{-1}\not\in R_{M}$, then $M\in {\mathcal{T}}_{t}(R)$, as
desired. \end{proof}

The following basic facts provide some background to the theorem and will be of use in its proof.
\begin{itemize}
\item {\tt Fact 1.} For each ideal $I$ of $R$, we have $I_{t} =
\bigcap_{M\in \Max_{t}(R)} IR_{M}$ \cite{Kg}, Theorem 2.9.

\item {\tt Fact 2.} Let $R$ be a Pr\"ufer domain, $I$ an ideal of $R$,
and $A$ and $B$  $R$-submodules of $\qf(R)$. Then $I(A\cap B)=IA\cap
IB$ \cite{Ba1}, Lemma 2.6.

\item {\tt Fact 3.} For a $t$-ideal $I$ of a domain $R$, let
$\overline{M}(R, I)=\{M\in \Max_{t}(R)| I\nsubseteq M\}$ and
$\mathcal{C}_{t}(I) =\bigcap_{M}R_{M}$ where $M$ ranges over
$\overline{M}(R, I)$. Then
$(\mathcal{C}_{t}(I):I)=\mathcal{C}_{t}(I)$. Indeed, it is clear
that $\mathcal{C}_{t}(I)\subseteq (\mathcal{C}_{t}(I):I)$.
Conversely, let $x\in (\mathcal{C}_{t}(I):I)$. For each $M\in
\overline{M}(R, I)$, let $a\in I\setminus M$. Since $xI\subseteq
\mathcal{C}_{t}(I)\subseteq R_{M}$, then $xa\in R_{M}$. So
$x=\frac{xa}{a}\in R_{M}$. Hence $x\in \mathcal{C}_{t}(I)$ and
therefore $(\mathcal{C}_{t}(I):I)=\mathcal{C}_{t}(I)$.

\item {\tt Fact 4.} For each $t$-ideal $I$ of a domain $R$ with finite
$t$-character, there exists a nonzero finitely generated subideal
$J$ of $I$ such that $\Max_{t}(R,I)=\Max_{t}(R,J)$. The proof apes
that of \cite{Ba1}, Lemma 2.13, by replacing ``maximal ideals" with
``$t$-maximal ideals."
\end{itemize}

\begin{proof}[Proof of Theorem~\ref{sec:3.2}] Assume $R$ is a PVMD
which is Clifford $t$-regular and let $0\not =x\in R$. We must show
that $\Max_{t}(R, x)$ is finite. Suppose by way of contradiction
that $\Max_{t}(R, x)$ is infinite. By Lemma~\ref{sec:3.9}, there is
$M\in \Max_{t}(R, x)\setminus {\mathcal{T}}_{t}(R)$. By
Lemma~\ref{sec:3.10}, $M=\bigcup_{N}M\wedge N$ where $N$ ranges over
$\Max_{t}(R)\setminus \{M\}$. Since $R_{M}$ is a valuation domain,
 $N$ may range over
$\Max_{t}(R,x)\setminus \{M\}$, so that $\{P_{\alpha}\}_{\alpha\in
\Omega}=\{M\wedge N\}_{M\not=N\in \Max_{t}(R,x)}$ is an infinite
totally ordered set. For each $\alpha\in \Omega$, we have
$0\subsetneq (x)\subseteq P_{\alpha}=M\wedge N_{\alpha}\subsetneq
N_{\alpha}$, for some $N_{\alpha}\in \Max_{t}(R,x)$. By
\cite{Kap}, Theorem 11, there exist distinct prime ideals
$P'_{\alpha}$ and $Q_{\alpha}$ such that $0\subsetneq
P_{\alpha}\subseteq P'_{\alpha}\subsetneq Q_{\alpha}\subseteq
N_{\alpha}$ with $\htt(Q_{\alpha}/P'_{\alpha})=1$.

\begin{claim}
For every $\alpha\not =\beta$, $Q_{\alpha}$ and $Q_{\beta}$ are
incomparable.
\end{claim}

We may assume $P_{\alpha}\subsetneq P_{\beta}$. Suppose that
$Q_{\alpha}\subseteq Q_{\beta}$. Then $Q_{\alpha}$ and $P_{\beta}$
are comparable. If $Q_{\alpha}\subseteq P_{\beta}$, then
$P_{\alpha}\subsetneq Q_{\alpha}\subseteq M\wedge
N_{\alpha}=P_{\alpha}$, absurd. If $P_{\beta}\subseteq Q_{\alpha}$,
then $P_{\beta}\subseteq M\wedge N_{\alpha}=P_{\alpha}\subsetneq
P_{\beta}$, absurd. Now, if $Q_{\beta}\subseteq Q_{\alpha}$, then
$P_{\beta}\subseteq M\wedge N_{\alpha}=P_{\alpha}$, which is absurd
too. This proves the claim.

Since $P_{\alpha}\subsetneq Q_{\alpha}$, then
$Q_{\alpha}\not\subseteq M$. For each $\alpha$, let $a_{\alpha}\in
Q_{\alpha}\setminus M$ and consider the ideal
$J_{\alpha}=P_{\alpha}+a_{\alpha}R$.

\begin{claim} $J_{\alpha}$ is
$t$-invertible.\end{claim}

By Lemma~\ref{sec:3.4}, it suffices to check that $J_{\alpha}$ is
$t$-locally principal. Let $N$ be a $t$-maximal ideal of $R$. Assume
-without loss of generality- that $J_{\alpha}\subseteq N$. Since
$R_{N}$ is a valuation domain and $a_{\alpha}\not\in P_{\alpha}$,
then $P_{\alpha}R_{N}\subsetneq a_{\alpha}R_{N}$. Hence
$J_{\alpha}R_{N}=a_{\alpha}R_{N}$, as desired. Therefore there is a
finitely generated subideal $F_{\alpha}$ of $J_{\alpha}$ such that
$(F_{\alpha})_{v}=(F_{\alpha})_{t}=(J_{\alpha})_{t}=(J_{\alpha})_{v}$.

Moreover, by Lemma~\ref{sec:3.7}, there is a finitely generated
subideal $I_{\alpha}$ of $Q_{\alpha}$ such that $\Max_{t}(R,
I_{\alpha})=\Max_{t}(R, Q_{\alpha})$. Consider the finitely
generated ideal given by $A_{\alpha}= F_{\alpha}+I_{\alpha}$. Since
$I_{\alpha}\subseteq A_{\alpha}\subseteq Q_{\alpha}$, then
$\Max_{t}(R, A_{\alpha})=\Max_{t}(R, Q_{\alpha})$. Finally, let
$B=\sum_{\alpha\in \Omega}(A_{\alpha})^{-1}$.

\begin{claim}  $B$ is a fractional ideal of $R$ which is $t$-invertible.\end{claim}

 Indeed, for each $\alpha$, we have $(x)\subseteq P_{\alpha}=(P_{\alpha})_{t}\subseteq
(J_{\alpha})_{t}=(F_{\alpha})_{t}\subseteq (A_{\alpha})_{t}$. So
$x(A_{\alpha})^{-1}\subseteq
(A_{\alpha})_{t}(A_{\alpha})^{-1}=(A_{\alpha})_{t}((A_{\alpha})_{t})^{-1}\subseteq
R$. Hence $xB\subseteq R$ and therefore $B$ is a fractional ideal of
$R$. Now let $N$ be a $t$-maximal ideal of $R$. {\bf Case 1:}
$A_{\alpha}\not\subseteq N$ for each $\alpha\in \Omega$. Since
$A_{\alpha}$ is finitely generated, then
$(A_{\alpha})^{-1}R_{N}=(A_{\alpha}R_{N})^{-1}=R_{N}$. Hence
$BR_{N}=R_{N}$. {\bf Case 2:} $A_{\alpha}\subseteq N$ for some
$\alpha\in \Omega$. Since $\Max_{t}(R, A_{\alpha})=\Max_{t}(R,
Q_{\alpha})$, then for each $\beta\not = \alpha$, $
A_{\beta}\not\subseteq N$. Otherwise, $N\in \Max_{t}(R, A_{\beta})=
\Max_{t}(R, Q_{\beta})$. Then $Q_{\alpha}$ and $Q_{\beta}$ are
comparable since both are included in $N$, absurd by the first
claim. Thus $N$ contains exactly one $A_{\alpha}$. So
$BR_{N}=(A_{\alpha})^{-1}R_{N}=(A_{\alpha}R_{N})^{-1}=a^{-1}R_{N}$
where $A_{\alpha}R_{N}=aR_{N}$ since $A_{\alpha}R_{N}$ is a finitely
generated ideal of the valuation domain $R_{N}$. It follows that $B$
is $t$-locally principal and therefore $t$-invertible
(Lemma~\ref{sec:3.4}).

Consequently, there is a finitely generated subideal $L$ of $B$ such
that $L_{v}=L_{t}= B_{t}=B_{v}$. There exist  ${\alpha}_{1}, \dots,
{\alpha}_{r}$ such that $L\subseteq \sum_{1\leq i\leq
r}(A_{\alpha_{i}})^{-1}\subseteq B$. Therefore $B^{-1}=L^{-1}=
(\sum_{1\leq i\leq r}(A_{\alpha_{i}})^{-1})^{-1}=\bigcap_{1\leq
i\leq r}(A_{\alpha_{i}})_{v} =\bigcap_{1\leq i\leq
r}(A_{\alpha_{i}})_{t}$. Now, let $\alpha\in \Omega\setminus
\{{\alpha}_{1}\dots, \alpha_{r}\}$. Then $\bigcap_{1\leq i\leq
r}(A_{\alpha_{i}})_{t}=B^{-1}\subseteq
(A_{\alpha})_{v}=(A_{\alpha})_{t}\subseteq Q_{\alpha}\subseteq
N_{\alpha}$. So there is $i\in \{1,\dots, r\}$ such that
$(A_{\alpha_{i}})_{t}\subseteq N_{\alpha}$. Hence $N_{\alpha}\in
\Max_{t}(R, A_{\alpha_{i}})= \Max_{t}(R, Q_{\alpha_{i}})$ and then
$Q_{\alpha_{i}}\subseteq N_{\alpha}$. This forces $Q_{\alpha}$ and
$Q_{\alpha_{i}}$ to be comparable, the desired contradiction. Thus
$\Max_{t}(R, x)$ is finite.

Next, we prove the converse of the theorem. Assume $R$ is a
Krull-type domain. Let $I$ be a $t$-ideal of $R$,
$\Max_{t}(R,I)=\{M_{1},\dots, M_{n}\}$ and $J=I^{2}(I:I^{2})$. We
wish to show that $I=J_{t}$. By Fact 1, it suffices to show that
 $IR_{M}=JR_{M}$ for each $t$-maximal ideal of $R$. Let $M\in
\Max_{t}(R)$. If $I\not\subseteq M$, then $J\not\subseteq M$ (since
$I^{2}\subseteq J$). So $IR_{M}= JR_{M}=R_{M}$. Assume  $I\subseteq
M$. Mutatis Mutandis, we may assume that $M=M_{1}$. One can easily
check via Fact 1 that $(I:I)=(\bigcap_{i=1}^{n}
(IR_{M_{i}}:IR_{M_{i}}))\cap \mathcal{C}_{t}(I)$. By Fact 3,
$(I:I^{2})=(\bigcap_{i=1}^{n} (IR_{M_{i}}:I^{2}R_{M_{i}}))\cap
\mathcal{C}_{t}(I)=(IR_{M_{1}}:I^{2}R_{M_{1}})\cap
(\bigcap_{i=2}^{n} (IR_{M_{i}}:I^{2}R_{M_{i}}))\cap
\mathcal{C}_{t}(I)$. Let $A=(IR_{M_{1}}:I^{2}R_{M_{1}})$ and
$B=(\bigcap_{i=2}^{n}(IR_{M_{i}}:I^{2}R_{M_{i}}))$. We have
$JR_{M_{1}}= I^{2}R_{M_{1}}(AR_{M_{1}}\cap BR_{M_{1}}\cap
\mathcal{C}_{t}(I)R_{M_{1}})$. By applying Fact 2 in the valuation
domain $R_{M_{1}}$, we obtain $JR_{M_{1}}=
(I^{2}R_{M_{1}}AR_{M_{1}})\cap (I^{2}R_{M_{1}}BR_{M_{1}})\cap
(I^{2}R_{M_{1}}\mathcal{C}_{t}(I)R_{M_{1}})$.

On one hand, $I^{2}R_{M_{1}}AR_{M_{1}}= IR_{M_{1}}$ since
$R_{M_{1}}$ is Clifford regular \cite{BS}. Further, we claim that
$I^{2}R_{M_{1}}BR_{M_{1}}\supseteq IR_{M_{1}}$. Indeed,
$I^{2}R_{M_{1}}BR_{M_{1}}=\bigcap_{i=2}^{n}I^{2}R_{M_{1}}(IR_{M_{i}}:I^{2}R_{M_{i}})R_{M_{1}}
=\bigcap_{i=2}^{n}\big(I^{2}(IR_{M_{i}}:I^{2}R_{M_{i}})\big)R_{M_{1}}
=\bigcap_{i=2}^{n}\big(I^{2}R_{M_{i}}(IR_{M_{i}}:I^{2}R_{M_{i}})\big)R_{M_{1}}
=\bigcap_{i=2}^{n}IR_{M_{i}}R_{M_{1}}\supseteq IR_{M_{1}}$, as
claimed; the first equality is due to Fact 2 and the last equality
holds because $R_{M_{i}}$ is Clifford regular.

On the other hand, $\mathcal{C}_{t}(I)R_{M_{1}}$ is an overring of
$R_{M_{1}}$ and hence $\mathcal{C}_{t}(I)R_{M_{1}}= R_{P}$ for some
$t$-prime ideal $P$ of $R$ contained in $M_{1}$. We claim that
$I\not\subseteq P$. Indeed, by Fact 4, there exists a nonzero
finitely generated ideal $L$ with $L\subseteq J_{t}\subseteq I$ and
$\Max_{t}(R,L)=\Max_{t}(R,J_{t})=\Max_{t}(R,I)$. So
$\mathcal{C}_{t}(I)=\mathcal{C}_{t}(J)=\mathcal{C}_{t}(L)$. Since
$R$ is integrally closed, $(L:L^{2})=((L:L):L)=(R:L)=L^{-1}$.
Furthermore it is easily seen that $L^{-1}\subseteq
\mathcal{C}_{t}(L)$. So $L^{2}(L:L^{2})\subseteq
L^{2}\mathcal{C}_{t}(L)=L^{2}\mathcal{C}_{t}(I)$. Since $R_{M_{1}}$
is Clifford regular, we get
$LR_{M_{1}}=L^{2}R_{M_{1}}(LR_{M_{1}}:L^{2}R_{M_{1}})=(L^{2}(L:L^{2}))R_{M_{1}}\subseteq
L^{2}\mathcal{C}_{t}(I)R_{M_{1}}=L^{2}R_{P}$. It results that
$LR_{P}\subseteq L^{2}R_{P}$ and hence $LR_{P}=L^{2}R_{P}$. By
\cite{Kap}, Theorem 76, $LR_{P}=R_{P}$. Hence $L\not\subseteq P$ and
thus $I\not\subseteq P$. This proves our claim

Now, using the above claims, we obtain
$JR_{M_{1}}=I^{2}R_{M_{1}}AR_{M_{1}}\cap
I^{2}R_{M_{1}}BR_{M_{1}}\cap
I^{2}R_{M_{1}}\mathcal{C}_{t}(I)R_{M_{1}}= IR_{M_{1}}\cap
I^{2}R_{P}=IR_{M_{1}}\cap R_{P}=IR_{M_{1}}$. Consequently,
$I=J_{t}$, as desired. This completes the proof of the theorem.
\end{proof}

Since in a Pr\"ufer domain the $t$-operation coincides with the
trivial operation, we recover Bazzoni's theorem (mentioned above) as
a consequence of Theorem~\ref{sec:3.2}. Recall at this point
Zanardo-Zannier's result
 that ``an integrally closed Clifford regular domain is
Pr\"ufer" \cite{ZZ}. Also it is worthwhile noticing that during the
proof of  Theorem~\ref{sec:3.2} we made use of Bazzoni-Salce result
that ``a valuation domain is Clifford regular" \cite{BS}.

\begin{corollary}[Bazzoni \cite{Ba4}, Theorem 4.5]\label{sec:3.11}
An integrally closed domain $R$ is Clifford regular if and only if
$R$ is a Pr\"ufer domain of finite character.
\end{corollary}

The next result solves Conjecture~\ref{sec:3.1} for the context of
strongly $t$-discrete domains.

\begin{corollary} \label{sec:3.12}
Assume $R$ is a strongly $t$-discrete domain. Then $R$ is a
pseudo-integrally closed Clifford $t$-regular domain if and only if
$R$ is a Krull-type domain.
\end{corollary}

\begin{proof} In view of Theorem~\ref{sec:3.2}, we only need to prove the ``only if"
assertion. Precisely, it remains to show that $R$ is a PVMD. Let $I$
be a finitely generated ideal of $R$. If $I_{t}=R$, then $I^{-1}=R$
and therefore $(II^{-1})_{t}=R$, as desired. Assume that $I_{t}$ is
a proper $t$-ideal of $R$. Suppose by way of contradiction that $I$
is not $t$-invertible. Let $M$ be a $t$-maximal ideal of $R$
containing $J=(II^{-1})_{t}$. Since $R$ is pseudo-integrally
closed, $(I_{t}:I_{t})=R$. Hence
$(I_{t}:(I_{t})^{2})=((I_{t}:I_{t}):I_{t})=(R:I_{t})=I^{-1}$.
Further $R$ is Clifford $t$-regular, then
$I_{t}=((I_{t})^{2}(I_{t}:(I_{t})^{2}))_{t}=(I_{t}^{2}I^{-1})_{t}=(IJ)_{t}$.
Therefore $R\subseteq J^{-1}=(J:J)\subseteq (IJ:IJ)\subseteq
((IJ)_{t}:(IJ)_{t})=(I_{t}:I_{t})=R$. Consequently,
$J^{-1}=(J:J)=R$. Hence $R\subseteq M^{-1}\subseteq J^{-1}=R$,
whence $M^{-1}=(M:M)=R$. So $(M:M^{2})=((M:M):M)=(R:M)=R$. Since $R$
is Clifford $t$-regular, then $M=(M^{2}(M:M^{2}))_{t}=(M^{2})_{t}$,
and hence $M$ is $t$-idempotent. This contrasts with the hypothesis
that $R$ is strongly $t$-discrete. It follows that $I$ is
$t$-invertible and thus $R$ is a PVMD. \end{proof}

\section{\bf Examples}\label{sec:4}

This section is motivated by an attempt to generating  original
families of integral domains with Clifford $t$-class
semigroup. Next, we announce our first result of this section. It
provides necessary and sufficient conditions for a pullback to
inherit the Krull type notion.

\begin{proposition}\label{sec:4.1}
Let $T$ be an integral domain, $M$ a maximal ideal of $T$, $K$ its
residue field, $\phi:T\longrightarrow K$ the canonical surjection,
and $D$ a proper subring of $K$. Let $R=\phi^{-1}(D)$ be the
pullback issued from the following diagram of canonical
homomorphisms:
\[\begin{array}{ccl}
R            & \longrightarrow                 & D\\
\downarrow   &                                 & \downarrow\\
T            & \stackrel{\phi}\longrightarrow  & K=T/M
\end{array}\]
Then $R$ is a Krull-type domain if and only if $D$ is a semilocal
B\'ezout domain with $\qf(D)=K$ and $T$ is a Krull-type domain such
that $T_{M}$ is a valuation domain.
\end{proposition}

\begin{proof}
By \cite{FG}, Theorem 4.1, $R$ is a PVMD if and only if so are $T$
and $D$, $\qf(D)=K$, and $T_{M}$ is a valuation domain. Now notice
that $T=S^{-1}R$, where $S=\phi^{-1}(D\setminus\{0\})$. Moreover,
by \cite{Kg}, Corollary 2.47, $P$ is a $t$-prime ideal of $R$ if and
only if $PT$ is a $t$-prime ideal of $T$, for every prime ideal $P$
of $R$ saturated with respect to $S$. Also,  by \cite{FG}, Proposition
1.8, $q$ is a $t$-maximal ideal of $D$ if and only if
$\phi^{-1}(q)$ is a $t$-maximal ideal of $T$, for every prime ideal
$q$ of $D$. Finally, if $A$ is a domain with only a finite number of
maximal $t$-ideals, then each maximal ideal of $A$ is a $t$-ideal
\cite{Z}, Proposition 3.5. Using the above four facts, we can easily
see that $R$ has finite $t$-character if and only if $D$ is a
semilocal B\'ezout domain and $T$ has finite $t$-character.\end{proof}

The next result investigates the transfer of the finite
$t$-character condition to polynomial rings.

\begin{proposition}\label{sec:4.2}
Let $R$ be an integrally closed domain and $X$ an indeterminate over
$R$. Then $R$ has finite $t$-character if and only if so does
$R[X]$.
\end{proposition}

\begin{proof} Assume that $R$ has finite $t$-character and let $f$
be a nonzero nonunit element of $R[X]$ and
$\{Q_{\alpha}\}_{\alpha\in\Omega}$ the set of all $t$-maximal ideals
of $R[X]$ containing $f$. Set $\Omega_{1}=\{\alpha\in\Omega\ |\
Q_{\alpha}\cap R=0\}$ and $\Omega_{2}=\{\alpha\in\Omega\ |\
q_{\alpha}=Q_{\alpha}\cap R\not=0\}$. Assume $\alpha\in\Omega_{1}$
and let $K=\qf(R)$ and $S=R\setminus\{0\}$. Then
$S^{-1}Q_{\alpha}$ is a maximal ideal of $K[X]$. Further $f$ is not
a unit in $K[X]$ since $Q_{\alpha}\cap R=0$. Now $K[X]$ is of finite
character (since a PID), then
$\{S^{-1}Q_{\alpha}\}_{\alpha\in\Omega_{1}}$ is finite (and so is
$\Omega_{1}$). Assume $\alpha\in\Omega_{2}$. By \cite{Kg}, Lemma
2.32, $q_{\alpha}$ is a $t$-prime ideal of $R$ with
$Q_{\alpha}=q_{\alpha}[X]$. We claim that $q_{\alpha}$ is
$t$-maximal in $R$. Deny. Then $q_{\alpha}\subsetneq M_{\alpha}$ for
some $M_{\alpha}\in\Max_{t}(R)$. So
$Q_{\alpha}=q_{\alpha}[X]\subsetneq M_{\alpha}[X]$, absurd since
$M_{\alpha}[X]$ is a $t$-prime ideal of $R[X]$. Now let $a$ denote
the leading coefficient of $f$. Clearly, $0\not= a\in q_{\alpha}$
(since $Q_{\alpha}=q_{\alpha}[X]$). Therefore
$\{q_{\alpha}\}_{\alpha\in\Omega_{2}}$ is finite (and so is
$\Omega_{2}$) since $R$ has finite $t$-character. Consequently,
$\Omega$ is finite, as desired. The converse lies on the fact that
the extension of a $t$-maximal ideal of $R$ is $t$-maximal in
$R[X]$.\end{proof}

Notice at this point that (as in Example~\ref{sec:2.7}) one can
build numerous examples of non-PVMD Clifford (or Boole) $t$-regular
domains through Propositions \ref{sec:2.4} or \ref{sec:2.5} combined
with Theorem~\ref{sec:2.6}. Next, we provide new families of
Clifford $t$-regular domains originating from the class
of PVMDs via a combination of Theorems \ref{sec:3.2} and Propositions \ref{sec:4.1} \& \ref{sec:4.2}.

\begin{example}\label{sec:4.3}
For each integer $n\geq 2$, there exists a
PVMD $R_{n}$ subject to the following conditions:\\
\1 $\dim(R_{n})=n$.\\
\2 $R_{n}$ is Clifford $t$-regular.\\
\3 $R_{n}$ is not Clifford regular.\\
\4 $R_{n}$ is not Krull.
\end{example}

\noindent Let $V_{0}$ be a rank-one valuation domain with
$K=\qf(V_{0})$. Let $V=K+N$ be a rank-one non strongly discrete
valuation domain (Cf. \cite{DP}, Remark 6(b)). We take
$R_{n}=V[X_{1}, \dots, X_{n-1}]$. For $n\geq 4$, the classical
$D+M$ construction provides more examples. Indeed, consider an
increasing sequence of valuation domains $V= V_{1}\subset
V_{2}\subset, \dots, \subset V_{n-2}$ such that, for each $i\in\{2,
\dots, n-2\}$, $\dim(V_{i})=i$ and $V_{i}/M_{i}=V/N= K$, where
$M_{i}$ denotes the maximal ideal of $V_{i}$. Set $T=V_{n-2}[X]$
and $M=(M_{n-2},X)$. Therefore $R_{n}=V_{0}+M$ is the desired
example. 


\begin{center}
\rule{4cm}{0.1mm}
\bigskip

\small
Department of Mathematical Sciences,
King Fahd University of Petroleum \& Minerals, P. O. Box 5046,
Dhahran 31261, Saudi Arabia\\
e-mail: kabbaj@kfupm.edu.sa\\
e-mail: amimouni@kfupm.edu.sa
\end{center}

\end{document}